\documentclass[12pt]{article}
\usepackage{amsmath,amssymb}

\newcommand{\Var}{\ensuremath{\mathcal{V}_{\mathbb{C}}}}

\def\uu{{\underline{u}}}

\def\kk{{\underline{k}}}

\def\1{\underline{1}}

\def\AA{\mathbb A}

\def\LLL{\mathbb L}
\def\Z{\mathbb Z}
\def\Q{\mathbb Q}

\def\H{\mathbb H}
\def\Ex{\mathbf{Exp}}
\def\Lo{\mathbf{Log}}

\newtheorem{theorem}{Theorem}

\newtheorem{proposition}{Proposition}

\newenvironment{definition}
{\smallskip\noindent{\bf Definition\/}:}{\smallskip\par}
\newenvironment{corollary}
{\smallskip\noindent{\bf Corollary\/}.}{\smallskip\par}

\newenvironment{remark}
{\smallskip\noindent{\bf Remark\/}.}{\smallskip\par}

\newenvironment{proof}
{\noindent{\bf Proof\/}.}{{ $\square$}\smallskip\par}

\title{Power structure over the Grothendieck ring of varieties and generating series of Hilbert schemes of points
\footnote{Math. Subject Class.: 14C05, 14G10}
}
\author{S.M.~Gusein-Zade \thanks{Partially supported by the grants
RFBR--04--01--00762,
NSh--1972.2003.1.
Address: Moscow State University,
Faculty of Mathematics and Mechanics, Moscow, 119992, Russia.
E-mail: sabir\symbol{'100}mccme.ru} \and
I.~Luengo
\and A.~Melle--Hern\'andez \thanks{The last two authors were partially
supported by the grant BFM2001--1488--C02--01. Address:
University Complutense de Madrid, Dept. of Algebra,
Madrid, 28040, Spain.
E-mail: iluengo\symbol{'100}mat.ucm.es, amelle\symbol{'100}mat.ucm.es}}
\date{}
\begin{document}
\def\eps{\varepsilon}

\maketitle

\begin{abstract}
The power structure over the Grothendieck (semi)ring of complex quasi-projective varieties constructed by the authors is used to express the generating series of classes of Hilbert schemes of zero-dimensional subschemes on a smooth quasi-projective variety of dimension $d$ as an exponent of that for the complex affine space $\AA^d$. Specializations of this relation give formulae for generating series of such invariants of the Hilbert schemes of points as the Euler characteristic and the Hodge--Deligne polynomial.
\end{abstract}

\section*{Introduction.}

The Grothendieck semiring $S_0(\Var)$ of complex quasi-projective varieties is the semigroup generated by isomorphism classes $[X]$ of such varieties modulo the relation $[X]=[X-Y]+[Y]$ for a Zariski closed subvariety $Y\subset X$; the multiplication is defined by the Cartesian product: $[X_1]\cdot [X_2]=[X_1\times X_2]$. The Grothendieck ring $K_0(\Var)$ is the group generated by these classes with the same relation and the same multiplication. Let $\LLL
\in K_0(\Var)$ be the class of the complex affine line, and let $K_0(\Var)[\LLL^{-1}]$ be the localization of Grothendieck ring $K_0(\Var)$ with respect to $\LLL$. A power structure over a (semi)ring $R$ (\cite{MRL}) is a map $\left(1+t\cdot R[[t]]\right)\times R\to 1+t\cdot R[[t]]$: $(A(t), m)\mapsto \left(A(t)\right)^m$ ($A(t)=1+a_1t+a_2t^2+\ldots$, $a_i\in R$, $m\in R$) such that all usual properties of the exponential function hold. In \cite{MRL} there was defined a power structure over each of the (semi)rings mentioned above. Here we use this structure to express the generating series of classes (in the Grothendieck (semi)ring of varieties) of Hilbert schemes of zero-dimensional subschemes on a smooth quasi-projective variety of dimension $d$ as an exponent of that for the complex affine space $\AA^d$. The conjecture that the generating series of Hilbert schemes of points on a smooth surface can be considered as an exponent was communicated to the authors by D.~van Straten. Specializations of this relation give formulae for generating series of certain invariants of the Hilbert schemes (Euler characteristic, Hodge--Deligne polynomial, ...) It is possible to define a power structure over the ring $\Z[u_1,\ldots, u_r]$ of polynomials in several variables with integer coefficients in such a way that, for $r=2$, it is the specialization of the power structure over the Grothendick ring $K_0(\Var)$ under the Hodge--Deligne polynomial homomorphism. This gives the main result of \cite{GS} and \cite{Cheah} (in somewhat different terms).

\section{Power structures.}
\begin{definition}
A {\em power structure} over a (semi)ring $R$ is a map
$\left(1+t\cdot R[[t]]\right)\times {R} \to 1+t\cdot R[[t]]$: $(A(t),m)\mapsto \left(A(t)\right)^{m}$,
which possesses the properties:
\begin{enumerate}
\item $\left(A(t)\right)^0=1$,
\item $\left(A(t)\right)^1=A(t)$,
\item $\left(A(t)\cdot B(t)\right)^{m}=\left(A(t)\right)^{m}\cdot
\left(B(t)\right)^{m}$,
\item $\left(A(t)\right)^{m+n}=\left(A(t)\right)^{m}\cdot
\left(A(t)\right)^{n}$,
\item $\left(A(t)\right)^{mn}=\left(\left(A(t)\right)^{n}\right)^{m}$,
\item $(1+t)^m=1+mt+$ terms of higher degree,
\item $\left(A(t^k)\right)^m = \left(A(t)\right)^m\raisebox{-0.5ex}{$\vert$}{}_{t\mapsto t^k}$.
\end{enumerate}
\end{definition}

\begin{remark}
In \cite{MRL} the properties~6 and 7 were not demanded, though the constructed power structures possessed them.
\end{remark}

\begin{definition}
A power structure is {\em finitely determined} if for each $N>0$ there exists $M>0$ such that the $N$-jet of the series $\left(A(t)\right)^m$ (i.e. $\left(A(t)\right)^m\, \mod t^{N+1}$) is determined by the $M$-jet of the series $A(t)$.
\end{definition}

One can see that one can take $M=N$.

We shall use the following general statement.

\begin{proposition}\label{prop1}
To define a finitely determined power structure over a ring $R$ it is sufficient to define the series $(1-t)^{-a}=1+at+\ldots$ for each $a\in R$ so that $(1-t)^{-(a+b)}=(1-t)^{-a}(1-t)^{-b}$, $(1-t)^{-1}=1+t+t^2+\ldots$
\end{proposition}

\begin{proof}
By properties 6 and 7 each series $A(t)\in 1+t\cdot R[[t]]$ can be in a unique way written as a product of the form $\prod\limits_{i=1}^\infty (1-t^i)^{-a_i}$ with $a_i\in R$. Then by properties 3 and 7 (and the finite determinacy of the power structure) one has 
\begin{equation}\label{eq1}
\left(A(t)\right)^m=\prod\limits_{i=1}^\infty (1-t^i)^{-a_i m}.
\end{equation}
In the other direction, one can easily see that the power structure defined by the equation~\ref{eq1} possesses the properties 1--7.
\end{proof}

Let $R_1$ and $R_2$ be rings with power structures over them. A ring homomorphism $\varphi:R_1\to R_2$ induces the natural homomorphism $R_1[[t]]\to R_2[[t]]$ (also denoted by $\varphi$) by $\varphi\left(\sum a_it^i\right)=\sum\varphi(a_i)t^i$. Proposition~\ref{prop1} yields the following statement.

\begin{proposition}\label{prop2}
If a ring homomorphism $\varphi:R_1\to R_2$ is such that $(1-t)^{-\varphi(a)}=\varphi\left((1-t)^{-a}\right)$, then
$\varphi\left(\left(A(t)\right)^m\right)=\left(\varphi\left(A(t)\right)\right)^{\varphi(m)}$.
\end{proposition}

One can describe the power structure over the Grothendieck ring $K_0(\Var)$ of complex quasi-projective varieties constructed in \cite{MRL} simply defining $(1-t)^{-[M]}$ for a quasi-projective variety $M$ as $\zeta_{[M]}(t)=1+[M]t+[S^2M]t^2+[S^3M]t^3+\ldots$, where $S^kM=M^k/S_k$ is the $k$th symmetric power of the variety $M$. The series $\zeta_{[M]}(t)$ is the Kapranov zeta function of the variety $M$:~\cite{Kap}. However, first, this does not permit to define the power structure over the Grothendieck semiring $S_0(\Var)$. (And one can say that elements of $S_0(\Var)$ have more geometric meaning: they are represented by "genuine" quasi-projective varieties, not by virtual one.) Second, the geometric description of the power structure has its own value (e.g. for Theorem~\ref{theo1} below). Moreover, one can say that the geometric construction of the power structure in \cite{MRL} gives a more fine definition of the coefficients of the series $\left(A(t)\right)^{[M]}$ preserving more structures on them. E.g., if the coefficients of the series $A(t)$ and the exponent $[M]$ are represented by compact spaces, coefficients of the series $\left(A(t)\right)^{[M]}$ also can be considered as such ones. It seems that the geometric construction of the power structure can be adapted for and used in some settings different from the Grothendieck ring of varieties. Therefore we shall describe the series $\left(A(t)\right)^{[M]}$, where $A(t)=1+\sum\limits_{i=1}^\infty[A_i]t^i$, $A_i$ and $M$ are quasi-projective varieties (in words a little bit different from those used in \cite{MRL}).

It is convenient to describe the power structure on the Grothendieck semiring $S_0(\Var)$ in terms of graded spaces (sets). A graded space (with grading from $\Z_{>0}$) is a space $A$ with a function $I_A$ on it with values in $\Z_{>0}$. The number $I_A(a)$ is called the weight of the point $a\in A$. To a series $A(t)=1+\sum\limits_{i=1}^\infty[A_i]\,t^i$ one associates the graded space $A=\coprod\limits_{i=1}^\infty A_i$ with the weight function $I_A$ which sends all points of $A_i$ to $i$. In the other direction, to a graded space $(A, I_A)$ there corresponds the series $A(t)=1+\sum\limits_{i=1}^\infty[A_i]\,t^i$ with $A_i=I_A^{-1}(i)$. To define the series $\left(A(t)\right)^{[M]}$, we shall describe the corresponding graded space $A^M$ first. The space $A^M$ consists of pairs $(K,\varphi)$, where $K$ is a finite subset of (the variety) $M$ and $\varphi$ is a map from $K$ to the graded space $A$. The weight function $I_{A^M}$ on $A^M$ is defined by $I_{A^M}(K,\varphi)=\sum\limits_{k\in K}I_A(\varphi(k))$. This gives a set-theoretic description of the series $\left(A(t)\right)^{[M]}$. To describe the coefficients of this series as elements of the Grothendieck semiring $S_0(\Var)$, one can write it as
$$
\left(A(t)\right)^{[M]}=1+
\sum_{k=1}^\infty
\left\{
\sum_{\kk:\sum ik_i=k}
\left[
\left(
(\prod_i M^{k_i})
\setminus\Delta
\right)
\times\prod_i A_i^{k_i}/\prod_iS_{k_i}
\right]
\right\}
\cdot t^k,
$$
where $\kk=\{k_i: i\in\Z_{>0}, k_i\in\Z_{\ge0}\}$, $\Delta$ is the "large diagonal" in $M^{\Sigma k_i}$ which consists of $(\sum k_i)$-tuples of points of $M$ with at least two coinciding ones, the permutation group $S_{k_i}$ acts by permuting corresponding $k_i$ factors in
$\prod\limits_i M^{k_i}\supset (\prod_i M^{k_i})\setminus\Delta$ and the spaces $A_i$ simultaneously (the connection between this formula and the description above is clear).

\section{Generating series of Hilbert schemes.}
Let $H^n_X$, $n\ge 1$, be the Hilbert scheme of zero-dimensional subschemes of length $n$ of a quasi-projective variety $X$; for $x\in X$, let $H^n_{X,x}$ be the Hilbert scheme of subschemes of $X$ concentrated at the point $x$. Let
$$
\H_X(t):=1+\sum\limits_{n=1}^\infty [H^n_X]\,t^n,\quad \H_{X,x}(t):=1+\sum\limits_{n=1}^\infty [H^n_{X,x}]\,t^n
$$
be the generating series of classes of Hilbert schemes $H^n_X$ and $H^n_{X,x}$ in the Grothendieck semiring $S_0(\Var)$. Let $\AA^d$ be the affine space of dimension $d$.

Computation of invariants of the Hilbert schemes $H^n_X$ for a smooth variety $X$ of dimension $n$ consists in two parts. First is computation of the corresponding invariants in the local case, i.e. invariants of the Hilbert schemes $H^n_{\AA^d,0}$ and the second is combining the local results to global ones. The following statements formalizes the second step (for invariants of classes in the Grothendieck (semi)ring). It generalizes the one for surfaces derived in \cite{MRL} from the computations of the motive of the Hilbert scheme of points on a surface by L.~G\"ottsche~\cite{Got2}.

\begin{theorem}\label{theo1}
For a smooth quasi-projective variety $X$ of dimension $d$,
\begin{equation}\label{eq2}
\H_X(t)=\left(\H^n_{\AA^d,0}(t)\right)^{[X]}.
\end{equation}
\end{theorem}

\begin{proof}
For a locally closed subvariety $Y\subset X$, let $H^n_{X,Y}$ be the Hilbert scheme of subschemes of $X$ concentrated at points of $Y$ and let $\H_{X,Y}(t):= 1+\sum\limits_{}^n [H^n_{X,Y}] t^n$ be the corresponding generating series. If $Y$ is a Zariski closed subset of $X$, then $\H_X(t)= \H_{X,Y}(t) \cdot \H_{X,X\setminus Y}(t)$. Therefore it is sufficient to prove that 
\begin{equation}\label{eq3}
\H_{X,Y}(t)=\left(\H^n_{\AA^d,0}(t)\right)^{[Y]}
\end{equation}
for a subvariety $Y$ of $X$ which lies in an affine space $\AA^N$ and such that the first $d$ affine coordinates of $\AA^N$ define local coordinates on $X$ at each point of $Y$. A zero-dimensional subscheme of $X$ concentrated at points of $Y$ is defined by a finite subset $K\subset Y$ to each point $x$ of which there corresponds a zero-dimensional subscheme of (the standard) affine space $\AA^d$ concentrated at the origin. The length of this subscheme is equal to the sum of lengths of the corresponding subschemes of $\AA^d$. Now (\ref{eq3}) follows immediately from the geometric description of the power structure over the Grothendieck semiring of quasi-projective varieties.
\end{proof}

Since $[\AA^d]=\LLL^d$, and therefore $\H^n_{\AA^d}(t)=\left(\H^n_{\AA^d,0}(t)\right)^{\LLL^d}$, one has the following statement.

\begin{corollary}
In $K_0(\Var)[\LLL^{-1}][[t]]$ one has:
$$
\H_X(t)=\left(\H^n_{\AA^d}(t)\right)^{\LLL^{-d}[X]}.
$$
\end{corollary}

Applying homomorphisms of power structures one can get specializations of the formula~(\ref{eq2}). The most known homomorphisms from the Grothendieck (semi)ring of quasi-projective varieties are the Euler characteristic $\chi$ (to the ring of integers $\Z$) and the Hodge--Deligne polynomial (to the ring $\Z[u,v]$ of polynomials in two variables). Over $\Z$ there is the standard power structure: the usual exponentiation. One has $\chi\left((1-t)^{-[X]}\right)=(1-t)^{-\chi(X)}$ (see, e.g., \cite{CDG}; this follows immediately from~\cite{Mac}), i.e., the Euler characteristic is a homomorphism of power structures. This implies the following statement.

\begin{proposition}\label{prop3}
For a smooth quasi-projective variety $X$ of dimension $d$,
$$
\chi\left(\H_X(t)\right)=\left(\chi(\H^n_{\AA^d,0}(t))\right)^{\chi(X)}.
$$
\end{proposition}

For $d=\dim X=2$ this (using \cite{ES}) gives 
$$
\chi\left(\H_X(t)\right)=\left(\prod\limits_{k\ge0}\frac{1}{1-t^k}\right)^{\chi(X)}.
$$
This is one of the results of \cite{Got1} obtained there using the Weil cojectures.

\section{Power structures on the ring of polynomials and Hodge--Deligne polynomials of Hilbert schemes.}
One can define a power structure over the ring $\Z[u_1,\ldots,u_r]$ of polynomials in $n$ variables with integer coefficients in the following way. Let $P(u_1,\ldots,u_r)=\sum\limits_{\kk\in\Z_{\ge0}^r}p_\kk \uu^\kk\in \Z[u_1,\ldots,u_r]$, where $\kk=(k_1,\ldots,k_r)$, $\uu=( u_1,\ldots,u_r)$, $\uu^\kk =u_1^{k_1}\cdot\ldots\cdot u_r^{k_r}$, $ p_\kk\in\Z$. Define $(1-t)^{-P(u_1,\ldots,u_r)}$ as $\prod\limits_\kk (1-\uu^\kk t)^{-p_\kk}$ where the power (with an integer exponent $-p_\kk$) means the usual one. One can easily see that $(1-t)^{-(P_1(\uu)+P_2(\uu))}= (1-t)^{-P_1(\uu)}(1-t)^{-P_2(\uu)}$ and therefore this defines a power structure over the ring $\Z[u_1,\ldots,u_r]$, i.e., for polynomials $A_i(\uu)$, $i\ge0$, and $M(\uu)$ there is defined a series $\left(1+A_1(\uu)t+A_2(\uu)t^2+\ldots\right)^{M(\uu)}$ with coefficients from $\Z[u_1,\ldots,u_r]$. Let $r=2$, $u_1=u$, $u_2=v$. Let $e: K_0(\Var)\to \Z[u,v]$ be the ring homomorphism which sends the class $[X]$ of a quasi-projective variety $X$ to its Hodge--Deligne polynomial $e_X(u,v)=\sum h_X^{ij}(-u)^i(-v)^j$. The following fact about Hodge--Deligne polynomials of symmetric powers of a variety is well known (see, e.g., \cite{Cheah}).

\begin{proposition}\label{prop4}
$$
e\left((1-t)^{-[X]}\right)=(1-t)^{-e_X(u,v)}
$$
where the powers are according to the power structures in the corresponding rings: $K_0(\Var)$ and $\Z[u,v]$ respectively.
\end{proposition}

Theorem~\ref{theo1} and Proposition~\ref{prop4} yield the following statement.

\begin{theorem}\label{theo2}
For a smooth quasi-projective variety $X$ of dimension $d$,
\begin{equation}\label{eq4}
e\left(\H_X(t)\right)=\left(e(\H^n_{\AA^d,0}(t))\right)^{e(X)}.
\end{equation}
\end{theorem}

This is the main result of \cite{GS} and \cite{Cheah} written in some sense in a more invariant way. To write it in the form similar to that of \cite{Cheah}, one should apply to the both sides of the equation~(\ref{eq4}) an isomorphism 
$L:1+t\cdot\Z[u,v][[t]]\to \Z[u,v][[t]]$ of abelian groups with multiplication and addition as group operations respectively. In \cite{Cheah} J.~Cheah used the usual logarithmic map $\log$ (and consequently the usual exponential map $\exp$ in the other direction). The same is done in a number of papers containing similar computations. These maps are defined only over the field $\Q$ of rational numbers what forces to write formulae in $\Q[u,v]$. Just in the same way one can use the (inverse to each other) isomorphisms $\Lo$ and $\Ex$ defined by 
$$
\Ex\left(P_1(\uu)t+P_2(\uu)t^2+\ldots\right):=\prod_{k\geq 1}(1-t^k)^{-P_k(\uu)}$$
(cf.~\cite{MRL}) or other ones: there are infinitely many of such pairs. This permits to write the formula staying in $\Z[u,v]$.

\end{document}